\documentclass[11pt]{amsart}
\usepackage[parfill]{parskip}    % Activate to begin paragraphs with an empty line rather than an indent
\usepackage{fullpage}
\usepackage{graphicx}
\usepackage[normalem]{ulem}
\usepackage{sseq}
\usepackage{amssymb,amsmath,latexsym}
\usepackage{tikz}
\usetikzlibrary{arrows,decorations.markings,decorations.pathreplacing,calc,matrix,intersections}
\tikzset{->-/.style={decoration={markings,mark=at position #1 with {\arrow{>}}},postaction={decorate}}}

\input epsf

\usepackage{epstopdf}
\usepackage{ifpdf}
\usepackage{color}
\usepackage{float}

\usepackage{dsfont}
\definecolor{red}{rgb}{1,0,0} 

 \definecolor{darkgreen}{rgb}{0, .7, 0}

 \definecolor{purple}{rgb}{.7, 0, 1}

\newcommand{\Z}{{\mathbb{Z}}}
\newcommand{\Q}{{\mathbb{Q}}}
\newcommand{\R}{{\mathbb{R}}}

\newcommand{\bdry}{\partial}
\newcommand{\inv}{^{-1}}
\newcommand{\pr}{^\prime}

\newtheorem{proposition}{Proposition}[section]
\newtheorem{definition}[proposition]{Definition}
\newtheorem{theorem}[proposition]{Theorem}
\newtheorem*{theorem*}{Theorem}
\newtheorem{lemma}[proposition]{Lemma}
\newtheorem*{claim}{Claim}
 
\newtheorem{corollary}[proposition]{Corollary}

\theoremstyle{remark}
\newtheorem{remark}[proposition]{Remark}

\title{The complex of free factors of a free group}
\author{Allen Hatcher$^1$ and Karen Vogtmann$^1$}
\thanks{$^1$Partially supported by NSF grant DMS-9307313} 
\begin{document}

\begin{abstract} We show that the geometric realization of the partially ordered set of proper free factors in a finitely generated free group of rank $n$ is homotopy equivalent to a wedge of spheres of dimension $n-2$.
\end{abstract}

\maketitle

{\it The original version of this paper, published in 1998 in the Oxford Quarterly, contained an error, and the purpose of the present version is to provide a correction. The main results of the paper are unchanged.  The error was in Lemma~2.3 which was used in the proof of Theorem~2.5.  Unfortunately, the statement of Lemma~2.3 was false;  it is replaced here with a new lemma which suffices to prove Theorem~2.5 after a few small changes in the proof.  Section~2 has been reorganized to accommodate the changes, and we have taken the occasion to make some further clarifications in this section.}

\section{Introduction}
An important tool in the study of the group $GL(n,\Z)$ is provided by the geometric realization of the partially ordered set (poset) of proper direct summands of $\Z^n$. The natural inclusion $\Z^n\to \Q^n$ gives a one-to-one correspondence between proper direct summands of $\Z^n$ and proper subspaces of $\Q^n$, so that this poset is isomorphic to the spherical building $X_n$ for $GL(n,\Q)$.  The term ``spherical" comes from the  Solomon-Tits theorem \cite{Sol}, which says that $X_n$ has the homotopy type of a bouquet of spheres:

\begin{theorem*} [Solomon-Tits Theorem]  {\sl The geometric realization of the poset of proper subspaces of an $n$-dimensional vector space has the homotopy type of a bouquet of spheres of dimension $n-2$. }
\end{theorem*}

The building $X_n$ encodes the structure of parabolic subgroups of $GL(n,\Q)$: they are the stabilizers of simplices.  $X_n$ also parametrizes the Borel-Serre boundary of the homogeneous space for $GL(n,\R)$.  The top-dimensional homology $H_{n-2}(X_n)$ is the Steinberg module $I_n$ for $GL(n,\Q)$, and is a dualizing module for the homology of $GL(n,\Z)$, i.e. for all coefficient modules $M$ there are isomorphisms
$$H^i(GL(n,\Z);M)\to H_{d-i}(GL(n,\Z);M\otimes I_n),$$
where $d=n(n-1)/2$ is the virtual cohomological dimension  of $GL(n,\Z)$.
 
If one replaces $GL(n,\Z)$ by the group $Aut(F_n)$ of automorphisms of the free group of rank $n$, the natural analog $FC_n$ of $X_n$ is the geometric realization of the poset of proper free factors of $F_n$. The abelianization map $F_n\to {\bf \Z}^n$ induces a map from $FC_n$ to the poset of summands of $\Z^n$.  In this paper we prove the analog of the Solomon-Tits theorem for $FC_n$:

\begin{theorem}   The geometric realization of the poset of proper free factors of $F_n$
has the homotopy type of a bouquet of spheres of dimension $n-2$.
\end{theorem}
 
By analogy, we call the top homology $H_{n-2}(FC_n)$ the {\it Steinberg module} for $Aut(F_n)$.  This leaves open some intriguing questions.  It has recently been shown that $Aut(F_n)$ is a virtual duality group~\cite{BF}; does the Steinberg module act as a dualizing module?  [This is answered in the negative for $n=5$ in the preprint \cite{HMNP} posted in 2022.]  There is an analog, called Autre space, of the homogeneous space for $GL(n,\Z)$ and the Borel-Serre boundary; what is the relation  between this and the ``building" of free factors?

In \cite{Qui}, Quillen developed tools for studying the homotopy type of the geometric realization $|X|$ of a poset $X$.  Given an order-preserving map $f\colon X\to Y$ (a ``poset map"), there is a spectral sequence relating the homology of $|X|$, the homology of $|Y|$, and the homology of the ``fibers" $|f/y|$,  where
$$f/y=\{x\in X|f(x)\leq y\}$$ 
with the induced poset structure. 

To understand $FC_n$ then, one might try to apply Quillen's theory using the poset map $FC_n\to X_n$.  However, it seems to be difficult to understand the fibers of this map.  Instead, we proceed by modeling the poset of free factors topologically, as the poset $B_n$ of simplices of a certain subcomplex of the ``sphere complex" $S(M)$ studied in \cite{Hat}. There is a natural poset map from $B_n$ to $FC_n$; we compute the homotopy type of $B_n$ and of the Quillen fibers of the  poset map, and apply Quillen's spectral sequence to obtain the result.

\section{Sphere systems}

Let $M$ be the compact 3-manifold obtained by taking a connected sum of $n$ copies of $S^1\times S^2$ and removing the interior of a closed ball.  A {\it sphere system} in $M$ is a non-empty finite set of disjointly embedded 2-spheres in the interior of $M$, no two of which are isotopic, and none of which bounds a ball or is isotopic to the boundary sphere of $M$.  The complex $S(M)$ of sphere systems in $M$ is defined to be the simplicial complex whose $k$-simplices are isotopy classes of sphere systems with $k+1$ spheres.  

Fix a basepoint $p$ on $\bdry M$. The fundamental group $\pi_1(M,p)$ is isomorphic to $F_n$.  Any automorphism of $F_n$ can be realized by a homeomorphism of $M$ fixing $\bdry M$.  A theorem of Laudenbach \cite{Lau} implies that such a homeomorphism inducing the identity on $\pi_1(M,p)$ acts trivially on isotopy classes of sphere systems, so that in fact $Aut(F_n)$ acts on $S(M)$.

For $H$ a subset of $\pi_1(M,p)$, define $S_H$ to be the subcomplex of $S(M)$ consisting of isotopy classes of sphere systems $S$ such that $\pi_1(M-S,p)\supseteq H$.  When $H$ is trivial, $S_H$ is $S(M)$, and in this case the following result was proved in \cite{Hat}. 

\begin{theorem}\label{SH} The complex $S_H$ is contractible for each $n\geq 1$.
\end{theorem}

The proof will be a variant of the proof in \cite{Hat}, using the following fact. 
 
\begin{lemma}\label{disjoint}  Any two simplices in $S_H$ can be represented by sphere systems $\Sigma$ and $S$ such that every element of $H$ is representable by a loop disjoint from both $\Sigma$ and $S$.
\end{lemma}

\begin{proof}  Enlarge $\Sigma$ to a maximal
sphere system $\Sigma'$, so the components of $M - \Sigma'$ are three-punctured spheres.  By Proposition~1.1 of \cite{Hat} we may isotope $S$ to be in normal form with respect to $\Sigma'$.  This means that $S$ intersects each component of $M - \Sigma'$ in a collection of surfaces, each having at most one boundary circle on each of the three punctures; and if one of these surfaces is a disk then it separates the two punctures not containing its boundary. 

We can represent a given element of $H$ by a loop $\gamma_0$ based at $p$, such that $\gamma_0$ is disjoint from $S$ and transverse to $\Sigma'$.  The points of intersection of $\gamma_0$ with $\Sigma\pr$ divide $\gamma_0$ into a finite set of arcs, each entirely contained in one component of $M-\Sigma\pr$.  Suppose one of these arcs $\alpha$, in a component $P$ of $M-\Sigma\pr$,  has both endpoints on the same boundary sphere $\sigma$ of $P$.  Since the map $\pi_0(\sigma - (S\cap \sigma)) \to \pi_0(P -(S\cap P))$ is injective (an easy consequence of normal form), there is an arc $\alpha'$ in $\sigma - (S\cap \sigma)$ with $\bdry \alpha' = \bdry\alpha$. Since $P$ is simply-connected, $\alpha$ is homotopic to $\alpha'$ fixing endpoints. This homotopy gives a homotopy of $\gamma_0$ eliminating the two points of $\bdry\alpha$ from $\gamma_0 \cap \Sigma'$, without introducing any intersection points with $S$.  After repeating this operation a finite number of times, we may assume there are no remaining arcs of $\gamma_0-(\gamma_0\cap\Sigma\pr)$ of the specified sort.

Now consider a homotopy $F \colon I \times I \to M$ of $\gamma_0$ to a loop $\gamma_1$ disjoint from $\Sigma$. Make $F$ transverse to $\Sigma\pr$ and look at $F\inv(\Sigma')$.  This consists of a collection of disjoint arcs and circles.  These do not meet the left and right edges of $I\times I$ since these edges map to the basepoint $p$.

We claim that every arc component of $F\inv(\Sigma\pr)$ with one endpoint on $I\times\{0\}$ must have its other endpoint on $I\times\{1\}$.  If not,  choose an ``edgemost" arc with both endpoints on $I\times \{0\}$, i.e. an arc such that the interval of $I\times\{0\}$ bounded by the endpoints contains no other point of $F\inv(\Sigma\pr)$.  Then $\gamma_0$ maps this interval to an arc $\alpha$ in $M-S$ which is entirely contained in one component $P$ of $M-\Sigma\pr$ and has both endpoints on the same boundary sphere of $P$, contradicting our
assumption that all such arcs have been eliminated.

Since the loop $\gamma_1$ is disjoint from $\Sigma$, it follows that $\gamma_0$ must be disjoint from $\Sigma$, and by construction $\gamma_0$ was disjoint from $S$.
\end{proof}

\begin{proof}[Proof of Theorem \ref{SH}]
Following the method in \cite{Hat}, a contraction of $S_H$ can be constructed by performing a sequence of surgeries on an arbitrary system $S$ in $S_H$ to eliminate its intersections with a fixed system $\Sigma$ in $S_H$, after first putting $S$ into normal form with respect to a maximal system $\Sigma'$ containing $\Sigma$.  In \cite{Hat} the system $\Sigma$ itself was maximal but for the present proof $\Sigma$ must be in $S_H$ so it cannot be maximal if $H$ is nontrivial.  We can choose $\Sigma$ to be a single sphere defining a vertex of $S_H$ for example.  Once $S$ has been surgered to be disjoint from $\Sigma$ it will lie in the star of $S$ which is contractible so this will finish the proof.  (Alternatively, we could use the simpler contraction technique of \cite{HaVo}, which reverses the roles of $S$ and $\Sigma$.) 
 
Each surgery on $S$ is obtained by taking a circle of intersection of $S$ and $\Sigma$ which is innermost on $\Sigma$ among the remaining circles of $S\cap \Sigma$, bounding a disk $D$ in $\Sigma$ with $D\cap S=\bdry D$, then taking the two spheres obtained by attaching parallel copies of $D$ to parallel copies of the two disks of $S-\bdry D$. By the lemma, elements of $H$ are representable by loops disjoint from $\Sigma$ and $S$, so these loops remain disjoint from sphere systems obtained by surgering $S$ along $\Sigma$ because such surgery produces spheres lying in a neighborhood of $S\cup\Sigma$. 

In order to ensure a continuous retraction the surgery process is made canonical by  performing surgery on all innermost spheres at once.  As a result, surgery can produce trivial spheres bounding  balls in $M$ or   balls punctured by the sphere $\bdry M$.   At the end of the surgery process we discard all trivial spheres in the resulting sphere system, and it must be checked that at least one nontrivial sphere remains. To check this we note that the end result could be achieved by doing a single surgery at a time and then renormalizing, so it suffices to  show that a single surgery cannot produce two trivial spheres. 

Suppose, to the contrary, that a nontrivial sphere $s$ is surgered to produced two trivial spheres $s'$ and $s''$. The spheres $s\cup s'\cup s''$ form the boundary of a three-punctured sphere $P$ in $M$.  If $s'$ or $s''$, say $s'$, bounds a ball or punctured ball on the same side of $s'$ as $P$ then $P$ will be contained in this ball or punctured ball, hence so will $s$, contradicting the nontriviality of $s$. Thus $s'$ and $s''$ both bound balls or punctured balls on the opposite side from $P$. They cannot both bound punctured balls since $M$ has only one puncture, so one of them bounds a ball. This forces the other one to be isotopic to $s$, but this is not possible since $s$ is nontrivial while $s'$ and $s''$ are trivial. 
\end{proof}

For an inductive argument in the next theorem we will need a generalization of the preceding theorem to manifolds with more than one boundary sphere. Let $M_k$ be the manifold  obtained from the connected sum of $n$ copies of $S^1 \times S^2$ by deleting the interiors of $k$ disjoint closed balls rather than just a single ball.  Choose the basepoint $p$ on one of the spheres in $\bdry M$. For $H\subseteq \pi_1(M_k,p)$, define $S_H(M_k)$ to be the complex of isotopy classes of sphere systems $S$ in $M_k$ no two of which are isotopic and none of which bounds a ball or is isotopic to a sphere of $\bdry M_k$, and such that $\pi_1(M-S,p)\supseteq H$. 
 
\begin{lemma}\label{SkH}
For $n\geq 1$ and $k\geq 1$ the complex $S_H(M_{k+1})$ deformation retracts onto a subcomplex isomorphic to $S_H(M_k)$.
\end{lemma}

\begin{proof}
When $H$ is trivial this is Lemma 2.2 in \cite{Hat} and the 
following proof extends the proof there to the general case after one small refinement to take $H$ into account.  Let the spheres of $\bdry M_{k+1}$ be $\bdry_0, \cdots,\bdry_k$ with $p\in \bdry_0$.  Call a vertex of $S_H(M_{k+1})$ {\it special\/} if it splits off a three-punctured sphere from $M_{k+1}$ having $\bdry_k$ as one of its boundary components.  Let $S'_H(M_{k+1})$ be the subcomplex of $S_H(M_{k+1})$ consisting of simplices with no special vertices. Then $S_H(M_{k+1})$ is obtained from $S'_H(M_{k+1})$ by attaching the stars of the special vertices to $S'_H(M_{k+1})$ along the links of the special vertices. These links can be identified with copies of $S_H(M_k)$ so they are contractible by induction on $k$.  The interiors of the stars are disjoint since there are no edges joining special vertices. Hence $S_H(M_{k+1})$ deformation retracts to $S'_H(M_{k+1})$ since it is obtained by attaching contractible complexes along contractible subcomplexes. 

The proof will be completed by showing that $S'_H(M_{k+1})$ deformation retracts onto a subcomplex isomorphic to $S_H(M_k)$.  Let $\Sigma$ be a sphere in $M_{k+1}$ splitting off a three-punctured sphere $P$ having $\bdry_0$ and $\bdry_k$ as its other two boundary components.  After putting a system $S$ in $S'_H(M_{k+1})$ into normal form with respect to a maximal system containing $\Sigma$, $S$ will intersect $P$ in a set of parallel disks separating $\bdry_0$ from $\bdry_k$.  We can eliminate these disks from $S\cap P$ by pushing them across $\bdry_k$ one by one and then outside $P$.  This can also be described as surgering $S$ along the circles of $S\cap \Sigma$ using the disks they bound in $\Sigma$ on the side of $\bdry_k$.  Such a surgery on a sphere $s$ of $S$ produces a pair of sphere $s'$ and $s''$ with $s'$ the sphere we have pushed across $\bdry_k$ and $s''$ a trivial sphere parallel to $\partial_k$ which is then discarded.

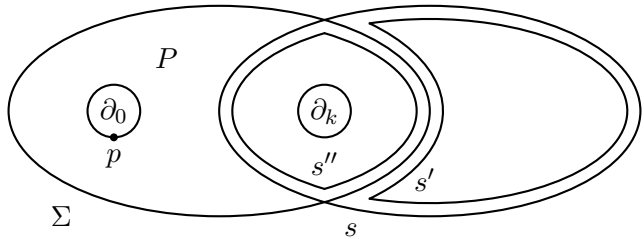
\begin{center}
  \begin{tikzpicture}[scale=.7] 
   \tikzset{partial ellipse/.style args={#1:#2:#3}{insert path={+ (#1:#3) arc (#1:#2:#3)} }}
 \draw [thick] (0,0) ellipse (4cm and 2cm); 
  \draw  [thick]   (4,0) ellipse (4cm and 2cm);  
    \draw [thick] (0,0) [partial ellipse=-58:58:3.75cm and 1.75cm];
          \draw [thick] (4,0) [partial ellipse=122:238:3.75cm and 1.75cm];
        \draw [thick] (0,0) [partial ellipse=-48:48:4.25cm and 2.25cm];
  \draw [thick] (4,0) [partial ellipse=-108:108:3.75cm and 1.75cm];
 \draw[thick] (-2,0) circle(.5);\node (b0) at (-2,0) {$\partial_0$};
\fill (-2,-.5) circle(.075); \node (p) at (-2,-.9 ) {$p$};
\draw[thick] (2,0) circle(.5);\node (b1) at ( 2,0) {$\partial_k$};
\node (sigma) at (-3,-2) {$\Sigma$};
\node (sigma) at (2.5,-2.25) {$s$};
\node (sigma) at ( 3.9,-1.3 ) {$s'$};
\node (sigma) at ( 2,-1) {$s''$};
\node (P) at (-1,1) {$P$};
  \end{tikzpicture} 
\end{center}

We claim that $s'$ is neither trivial nor special.  
If it were either of these, it would separate $M_{k+1}$ into two components, hence $s$ would also separate. The effect of replacing $s$ by $s'$ is to move the puncture $\bdry_k$ from one side of $s$ to the other side. It is then not hard to check that $s$ being neither trivial nor special implies the same is true for $s'$.  

The surgery defines a path in $S(M_{k+1})$ from $S$ to $S\cup s'$ and then to $(S\cup s')-s$.  The fundamental group of the component of $M_{k+1}-S$ containing $p$ is unchanged during this process so the path lies in $S'_H(M_{k+1})$.  By eliminating all the disks of $S\cap P$ by surgeries in this way we see  that $S'_H(M_{k+1})$ deformation retracts onto its subcomplex of the systems disjoint from $\Sigma$. This subcomplex can be identified with $S_H(M_k)$ by identifying $M_{k+1}-P$ with $M_k$ and choosing a new basepoint $p$ in $\Sigma$.  
\end{proof}

For the proof of the main theorem in the paper we will work with certain subcomplexes of $S(M_k)$ and $S_H(M_k)$,  the subcomplexes  $Y(M_k)\subset S(M_k)$ and $Y_H(M_k)\subset S_H(M_k)$ consisting of sphere systems $S$ with $M_k - S$ connected. Eventually only the case $k=1$ will be needed, but to prove the key property that $Y(M_1)$ and $Y_H(M_1)$ are highly connected we will need to consider larger values of $k$.  It will be convenient to extend the definition of $Y(M_k)$ and $Y_H(M_k)$ to allow $n=0$, with $M_k$ the sphere $S^3$ with $k$ punctures.  In this case $Y(M_k)$ and $Y_H(M_k)$ are empty since all spheres in $M_k$ are separating.

\begin{definition}  A simplicial complex $K$ is {\it $m$-spherical} if it is $m$-dimensional and $(m-1)$-connected.   A complex is {\it spherical\/} if it is $m$-spherical for some $m$.
\end{definition}

\begin{theorem}\label{YH} Let $H$ be a free factor of $F_n=\pi_1(M_k,p)$.  Then
$Y_H(M_k)$ is $(n-rk(H)-1)$-spherical, where $rk(H)$ is the rank of $H$.
\end{theorem}

In particular, when $H$ is trivial $Y(M_k)$ is $(n-1)$-spherical.  This special case is part of Proposition~3.1 of \cite{Hat} whose proof contained the same error that is corrected below.  A corrected proof of the special case already appeared in Proposition~3.2 of \cite{HW}.

\begin{proof}   
$Y_H(M_k)$ has dimension $n-rk(H)-1$ since a maximal simplex of $Y_H(M_k)$ has $n-rk(H)$ spheres. This follows because all free factors of $F_n$ of the same rank are equivalent under automorphisms of $F_n$, and all sphere systems with connected complement and the same number of spheres are equivalent under orientation-preserving homeomorphisms of $M_k$. Thus it suffices to prove $Y_H(M_k)$ is $(n-rk(H)-2)$-connected.  

Let $i\leq n-rk(H)-2$.  Any map  $g\colon S^i\to Y_H(M_k)$ can be extended to a map $\widehat g\colon D^{i+1}\to S_H(M_k)$ since $S_H(M_k)$ is contractible.  We can assume $\widehat g$ is a simplicial map with respect to some triangulation of $D^{i+1}$ compatible with its standard piecewise linear structure.  We will repeatedly redefine $\widehat g$ on the stars of certain simplices in the interior of $D^{i+1}$ until eventually the image of $\widehat g$ lies in $Y_H(M_k)$.  

To each sphere system $S$ we associate a dual graph $\Gamma(S)$, with one vertex for each component of $M-S$ and one edge for each sphere in $S$.  The endpoints of the edge corresponding to $s\in S$ are the vertices corresponding to the component or components adjacent to $s$. We say a sphere system  $S$ is   {\it purely separating\/} if $\Gamma(S)$ has no edges which begin and end at the same vertex.  Each sphere system $S$ has a {\it purely separating core}, consisting of those spheres in $S$ which correspond to the core of $\Gamma(S)$, i.e., the subgraph  spanned by edges with distinct vertices.  The purely separating core of $S\in S_H(M_k)$ is empty if and only if $S$ is in $Y_H(M_k)$.

Let $\sigma$ be a simplex of $D^{i+1}$ of maximal dimension among the simplices $\tau$ with $\widehat g(\tau)$ purely separating. Note that all such simplices $\tau$ lie in the interior of $D^{i+1}$ since the boundary of $D^{i+1}$ maps to $Y_H$.  
Let $S = \widehat g(\sigma)$, and let $N_0,\cdots,N_r$ ($r\geq 1$) be the connected components of $M-S$, with $p\in N_0$.  A simplex $\tau$ in the link $lk(\sigma)$  maps to a system  $T$ in the link of $S$, so that each $T_j=T\cap N_j$ is a sphere system in $N_j$ and $H\leq \pi_1(N_0-T_0,p)$.  Furthermore $N_j-T_j$ must be connected for all $j$ since otherwise the core of $\Gamma(S\cup T)$ would have more edges than $\Gamma(S)$, contradicting the maximality of $\sigma$.  Thus $\widehat g$ maps  $lk(\sigma)$ into a subcomplex of $S_H(M_k)$ which can be identified with $Y_H(N_0)*Y(N_1)*\cdots *Y(N_r)$.  Some of the factors $Y_H(N_0)$ and $Y(N_j)$ can be empty if $rk(H)=rk(\pi_1(N_0))$ or $rk(\pi_1(N_j))=0$. Such factors are $(-1)$-spherical and contribute nothing to the join. 

Since $\sigma$ is a simplex in the interior of $D^{i+1}$,   $lk(\sigma)$ is a sphere of dimension $i-dim(\sigma)$.  Each $N_j$ has fundamental group of rank $n_j\leq n$ with equality only if $N_j$ has fewer than $k$ boundary components, so by  
induction on the lexicographically ordered pair $(n,k)$, $Y_H(N_0)$ is $(n_0 - rk(H) - 1)$-spherical and, for $j\geq 1$, $Y(N_j)$ is $(n_j-1)$-spherical.  The induction can start with the cases $(n,k)=(0,k)$ when the theorem is obvious.  
For the join $Y_H(N_0)*Y(N_1)*\cdots *Y(N_r)$ it then follows that this is spherical of dimension $(\sum_{i=0}^r n_j) -rk(H)-1$.

Now  $  n=(\sum_j n_j) + rk(\pi_1(\Gamma(S))) = (\sum_j n_j) + m - r$ where $m$ is the number of spheres in $S$, i.e., edges in $\Gamma(S)$.  Since a simplicial map cannot increase dimension, we have $\dim(\sigma) \ge m-1$. Therefore  

\vskip-22pt

\begin{align*}
i-\dim(\sigma)&\leq n-rk(H)-2-\dim(\sigma) \\
        &\le n-rk(H)-m-1 \\
         &=\Bigl(\sum_j n_j\Bigr)-rk(H)-1-r \\ 
         & <\Bigl(\sum_j n_j\Bigr)-rk(H)-1. 
 \end{align*}
 
 \vskip-8pt
 
Hence the map $\widehat g\colon lk(\sigma)\to Y_H(N_0)*Y(N_1)*\cdots *Y(N_r)$ can be extended to a map of a disk $D^k$  into $Y_H(N_0)*Y(N_1)*\cdots *Y(N_r)$, where $k=i+1-dim(\sigma)$. The system $S$ is compatible with every system in the image of  $D^k$, so this map can be extended to a map $\sigma*D^k\to S_H$.   We replace the star of $\sigma$ in $D^{i+1}$ by the disk $\bdry(\sigma)*D^k$, and define $\widehat g$ on $\bdry(\sigma)*D^k$ using this map.  

What have we improved?  The new simplices in the disk $\bdry(\sigma)*D^k$ are of the form $\sigma\pr*\tau$, where $\sigma\pr$ is a face of $\sigma$ and $\widehat g(\tau)\subset Y_H(N_0)*Y(N_1)*\cdots *Y(N_r)$.  The image of such a simplex $\sigma\pr*\tau$ is a system $S\pr\cup T$ such that in $\Gamma(S\pr\cup T)$ the edges corresponding to $T$ are all loops. Therefore any simplex in the disk $\bdry(\sigma)*D^k$ with purely separating image must lie in the boundary of this disk, where we have not modified $\widehat g$. 

We continue this process, eliminating purely separating simplices until there are none in the image of $\widehat g$.  Since every system in $S_H(M_k)-Y_H(M_k)$ has a non-trivial purely separating core, in fact the whole disk maps into $Y_H(M_k)$, and we are done.  
\end{proof}

\section{Free factors}

We now turn to the poset  $FC_n$  of proper
free factors of the free group $F_n$, partially ordered by inclusion.   A $k$-simplex in the geometric realization $|FC_n|$ is a flag  $H_0<H_1<\cdots< H_k$ of proper free factors of $F_n$, each properly included in the next.  Each $H_i$ is also a free factor of $H_{i+1}$ (see ~\cite{MKS}, p. 117]), so that a maximal simplex of $|FC_n|$ has dimension
$n-2$.

We want to model free factors of $F_n$ by sphere systems in $Y=Y(M)$, by taking the fundamental group of the (connected) complement. Here $M$ is the manifold obtained from the connected sum of $n$ copies of $S^1\times S^2$ by removing the interior of a closed ball. A sphere system with $n$ spheres and connected complement, corresponding to an $(n-1)$-dimensional simplex of $Y$, in fact has simply-connected
complement.  But we only want to consider proper free factors, so instead we consider the $(n-2)$-skeleton $Y^{(n-2)}$.  Since $Y$ is $(n-2)$-connected by Theorem~\ref{YH}, $Y^{(n-2)}$ is $(n-2)$-spherical. 

In order to relate $Y^{(n-2)}$ to $FC_n$,  we take the barycentric subdivision $B_n$ of $Y^{(n-2)}$.  Then $B_n$ is the geometric realization of a poset of isotopy classes of sphere systems, partially ordered by inclusion. If $S\subseteq S\pr$ are sphere systems, we have
$\pi_1(M-S,p)\geq \pi_1(M-S\pr,p)$, reversing the partial ordering.  Taking fundamental group of the complement thus gives a poset map $f\colon B_n\to (FC_n)^{op}$, where $(FC_n)^{op}$ denotes $FC_n$ with the opposite partial ordering.

\begin{proposition}\label{lift}   $f\colon B_n\to (FC_n)^{op}$ is
surjective.
\end{proposition}

\begin{proof}  Every simplex of $FC_n$ is contained in a simplex of dimension $n-2$ so it suffices to show $f$ maps onto all $(n-2)$-simplices. The group $Aut(F_n)$ acts transitively on $(n-2)$-simplices of $FC_n$, and all elements of $Aut(F_n)$ are realized by homeomorphisms of $M$, so $f$ will be surjective if its image contains a single $(n-2)$-simplex, which it obviously does.
\end{proof}

\begin{corollary}\label{connected} $FC_n$ is connected if $n\geq 3$.
\end{corollary}

\begin{proof}   Theorem~\ref{YH} implies that $B_n$ is connected for $n\geq 3$.  So, given any two vertices of $FC_n$, lift them to  vertices of $B_n$ by Proposition~\ref{lift},  connect the lifted vertices by a path, then project the path back down to $FC_n$. 
\end{proof}  

For any proper free factor $H$, let  $B_{\geq H}$ denote the fiber $f/H$, consisting of isotopy classes of sphere systems $S$ in $B_n$ with $\pi_1(M-S,p)\geq H$.

\begin{proposition}\label{fiber} 
Let $H$ be a proper free factor of $F_n$.  Then $B_{\geq H}$ is $(n-rk(H)-1)$-spherical.  If $rk(H)=n-1$ then $B_{\geq H}$ is a single point.
\end{proposition}

\begin{proof} The fiber
$B_{\geq H}$ is the barycentric subdivision of $Y_H$, so is $(n-rk(H)-1)$-spherical by Theorem~\ref{YH}.  

Suppose $rk(H)=n-1$, so that $\pi_1(M,p)=H*\langle x\rangle$ for some $x$. An element of $B_{\geq H}$ is represented by a sphere system containing exactly one sphere $s$, which is non-separating with $\pi_1(M-s,p)=H$.  Suppose $s$ and $s\pr$ are two such spheres.  Since $s$ and $s\pr$ are both non-separating,  there is a homeomorphism $h$ of $M$ taking $s$ to $s\pr$. Since automorphisms of $H$ can be realized by homeomorphisms of $M$ fixing $s$, we may assume that the induced map on $\pi_1(M,p)$ is the identity on $H$.

\begin{claim}  The induced map $h_*\colon \pi_1(M,p)\to\pi_1(M,p)$ must send $x$ to an element of
the form
$Ux^{\pm 1} V$, with $U,V\in H$.
\end{claim}

\begin{proof}  Let $\{x_1,\cdots,x_{n-1}\}$ be a basis for $H$, and let $W$ be the reduced word representing $h_*(x)$ in the basis $\{x_1,\cdots,x_{n-1},x\}$ for
$\pi_1(M,p)$.  By looking at the map induced by $h$ on homology, we see that the exponent sum of $x$ in $W$ must be $\pm 1$.  Since $h_*$ fixes $H$, $\{x_1,\cdots,x_{n-1},W\}$ is a basis for $\pi_1(M,p)$.  If $W$ contained both $x$ and $x\inv$, we could apply Nielsen automorphisms to the set $\{x_1,\cdots,x_{n-1},W\}$ until $W$ was of the form $x^{\pm 1}W_0x^{\pm 1}$.    But $\{x_1,\cdots,x_{n-1},x^{\pm 1}W_0x^{\pm 1}\}$ is not a basis, since it is Nielsen reduced and not of the form $\{x_1^{\pm 1},\cdots,x_{n-1}^{\pm 1},x^{\pm 1}\}$ (see \cite{LySc}, Prop. 2.8)). 
\end{proof}

The automorphism fixing $H$ and sending $x\mapsto Ux^{\pm 1} V$ can be realized by a homeomorphism $h\pr$ of $M$ which takes $s$ to itself (see \cite{Lau}, Lemme 4.3.1).  The composition $h\pr h\inv$ sends $s\pr$ to $s$ and induces the identity on $\pi_1$, hence acts trivially on the sphere complex.  Thus $s$ and $s\pr$ are isotopic. 
\end{proof}

\begin{corollary}\label{oneconn}  $FC_n$ is simply connected for $n\geq 4$.
\end{corollary}

\begin{proof}  Let $e_0,e_1,\cdots,e_k$ be the edges of an edge-path loop in $FC_n$, and choose lifts $\tilde e_i$ of these edges to $B_n$. Let $e_{i-1}e_{i}$ or $e_ke_0$ be two adjacent edges of the path, meeting at the vertex $H$.  The lifts $\tilde e_{i-1}$ and $\tilde e_{i}$ may not be connected, i.e. $\tilde e_{i-1 }$ may terminate at a sphere system $S\pr$ and $\tilde e_i$ may begin at a different sphere system $S$.  However, both $S$  and $S\pr$ are in the fiber $B_{\geq H}$, which is connected by Proposition~\ref{fiber}, so we may  connect $S$ and $S\pr$ by a path in $B_{\geq H}$.  Connecting the endpoints of each lifted edge in this way, we obtain a loop in $B_n$, which may be filled in by a disk if $n\geq 4$, by Proposition~\ref{fiber}.   The projection of this loop to $FC_n$ is homotopic to the original loop, since each extra edge-path segment we added projects to a loop in the star of some vertex $H$, which is contractible.  Therefore the projection of the disk kills our original loop in the fundamental
group.
\end{proof}

\begin{remark} It is possible to describe the complex $Y_H$ purely in terms of $F_n$.  Suppose first that $H$ is trivial. Define a simplicial complex $Z$ to have vertices the rank $n-1$ free factors of $F_n$, with a set of $k$ such factors spanning a simplex in $Z$ if there is an automorphism of $F_n$ taking them to the $k$ factors obtained by deleting
the standard basis elements $x_1, \cdots, x_k$ of $F_n$ one at a time. There is a simplicial map $f\colon Y \to Z$ sending a system of $k$ spheres to the set of $k$ fundamental groups of the complements of these spheres. These fundamental groups are equivalent to the standard set of $k$ rank $n-1$ factors under an automorphism of $F_n$ since the homeomorphism group of $M$ acts transitively on simplices of $Y$ of a given
dimension, and the standard $k$ factors are the fundamental groups of the complements of the spheres in a standard system in $Y$.  The last statement of Proposition~\ref{fiber} says that $f$ is a bijection on vertices, so $f$ embeds $Y$ as a subcomplex of $Z$. The maximal simplices in $Y$ and $Z$ have dimension $n-1$ and the groups $Homeo(M)$ and
$Aut(F_n)$ act transitively on these simplices, so $f$ must be surjective, hence an isomorphism. When $H$ is nontrivial, $f$ restricts to an isomorphism from $Y_H$ onto the
subcomplex $Z_H$ spanned by the vertices which are free factors containing $H$.
\end{remark}

We are now ready to apply Quillen's spectral sequence to compute the homology of $B_n$ and thus prove the main theorem.

\begin{theorem}\label{spherical}  $FC_n$ is $(n-2)$-spherical.
\end{theorem}

\begin{proof}  We prove the theorem by induction on $n$.  If $n\leq 4$, the theorem follows from Corollaries~\ref{connected} and \ref{oneconn}.  

Quillen's spectral sequence (\cite{Qui}, 7.7)  applied to $f\colon B_n\to (FC_n)^{op}$ becomes 
$$E^2_{p,q}=H_p(FC_n;H\mapsto H_q(B_{\geq H}))\Rightarrow H_{p+q}(B_n),$$ 
where the $E^2$-term is computed using homology with coefficients in the functor $H\mapsto H_q(B_{\geq H})$. 

For $q=0$, Corollary~\ref{connected} gives $H_0(B_{\geq H})= \Z$ for all $H$, so $E^2_{p,0}=H_p(FC_n,\Z)$.

For $q>0$, we have $E^2_{p,q}=H_p(FC_n;H\mapsto \widetilde H_q(B_{\geq H}))$, and we follow Quillen (\cite{Qui}, proof of Theorem 9.1) to compute this.

For a subposet $A$ of $FC_n$, let $L_A$ denote the functor sending $H$ to a fixed abelian group $L$ if $H\in A$ and to $0$ otherwise.  Set $U=FC_{\leq H}=\{H\pr\in FC_n| H\pr\leq H\}$ and $V=FC_{<H}=\{H\pr\in FC_n| H\pr< H\}$.  Then 
\begin{align*}
 H_i(FC_n;L_V)&=H_i(V;L), \hbox{ and}\\
H_i(FC_n;L_U)&=H_i(U;L)=\begin{cases}L &\hbox{ if } $i=0$ \hbox{ and }\\ 0 &\hbox {otherwise},\end{cases}
\end{align*} 
since $|U|$ is contractible.  The short exact sequence of functors
$$
1\to L_V\to L_U\to L_{\{H\}}\to 1
$$
gives a long exact homology sequence, from which we compute
$$
H_i(FC_n;L_{\{H\}})=\widetilde H_{i-1}(FC_{<H};L).
$$
Now, $H\mapsto \widetilde H_q(B_{\geq H})$ is equal to the functor
$$
\bigoplus_{rk(H)=n-q-1} \widetilde H_q(B_{\geq H})_{\{H\}}
$$
since $B_{\geq H}$ is $(n-rk(H)-1)$-spherical by Proposition~\ref{fiber}.  Thus 
\begin{align*}
E^2_{p,q}&=H_p(FC_n;H\mapsto\widetilde H_q(B_{\geq H}))\\
                    &=\bigoplus_{rk(H)=n-q-1} H_p(FC_n;\widetilde H_q(B_{\geq H})_{\{H\}})\\
                    &=\bigoplus_{rk(H)=n-q-1} \widetilde H_{p-1}(FC_{<H};\widetilde H_q(B_{\geq H}))
\end{align*}

Free factors of $F_n$ contained in $H$ are also free factors of $H$.  Since $H$ has rank $<n$.  $FC_{<H}$ is $(rk(H)-2)$-spherical by induction.  Therefore $E^2_{p,q}=0$ unless  $p-1=(n-q-1)-2$, i.e. $p+q=n-2.$  Since all terms in the $E^2$-term of the spectral sequence are zero except the bottom row for $p\leq n-2$ and the diagonal $p+q=n-2$, all differentials are zero and we have $E^2=E^\infty$ as in the following diagram:

\begin{center}\begin{tikzpicture}
\matrix (m) [matrix of math nodes,
             nodes in empty cells,
             nodes={minimum width=2ex, minimum height=2ex,
                    text depth=1ex,
                    inner sep=0pt, outer sep=0pt,
                    anchor=base},
             column sep=1.5ex, row sep=1ex]%
{
n-2 & E^2_{0,n-2} &0 &0 & 0 & \cdots \\
n-3 &0& E^2_{1,n-3} & 0 & 0 & \cdots \\
n-4 &0& 0 & E^2_{2,n-4} & 0 & \cdots \\
\vdots    & \vdots       & \vdots       & \vdots & \ddots       &0\\
1    & 0 & 0& 0  &  \cdots &E^2_{n-3,1}& 0\\
0    & H_0(FC_n) & H_1(FC_n) & H_2(FC_n)  &  \cdots &H_{n-3}(FC_n)&H_{n-2}(FC_n)& 0\\
      &  0     &  1          &  2           & \cdots & n-3&n-2          & \strut \\
};
%\draw[-stealth, densely dotted] (m-3-4) to (m-2-2.east);  %d^2 differential. 
\draw[thick] (-6,-1.65) to (6 ,-1.65); %bottom border
 \draw[thick] (-5,2.5) to (-5,-2.25); %left border
\end{tikzpicture}
\end{center}

\vskip-6pt

Since $FC_n$ is connected and the spectral sequence converges to $H_*(B_n)$, which is $(n-3)$-connected, we
must have $\tilde H_i(FC_n)=0$ for $i\neq n-2$.  Since $FC_n$ is simply-connected by Corollary~\ref{oneconn}, this implies that $FC_n$ is $(n-3)$-connected by the Hurewicz theorem. The theorem follows since $FC_n$ is $(n-2)$-dimensional.\end{proof}

\section{The Cohen-Macaulay Property}

In a PL triangulation of an $n$-dimensional sphere, the link of every $k$-simplex is an $(n-k-1)$-sphere. A poset is said to be {\it Cohen-Macaulay} of dimension $n$ if its geometric realization is $n$-spherical and the link of every $k$-simplex is
$(n-k-1)$-spherical (see \cite{Qui}).   Spherical buildings are Cohen-Macaulay, and we remark that $FC_n$ also has this nice local property.

To see this, let $\sigma=\{H_0<H_1<\cdots<H_k\}$ be a $k$-simplex of $FC_n$.  The link of $\sigma$ is the join of subcomplexes  $FC_{H_i,H_{i+1}}$ of $FC_n $ spanned by free factors $H$ with $H_i<H<H_{i+1}$ ($-1\leq i\leq k$, with the conventions $H_{-1}=1$ and $H_{k+1}=F_n$).  Counting dimensions, we see that it suffices to show that for each $r$ and $s$ with $0\leq s <r$, the poset $FC_{r,s}$ of proper free factors $H$ of $F_r$ which properly contain $F_s$ is $(r-s-2)$-spherical.  The proof of this is identical to the proof that $FC_n$ is $(n-2)$-spherical, after setting $n=r$ and replacing the complex $Y$ by $Y_{F_s}$.

\section{The map to the building}

As mentioned in the introduction,
the abelianization map $F_n\to \Z^n$ induces a map from the free factor complex $FC_n$ to the building $X_n$, since summands of $\Z^n$ correspond to subspaces of $\Q^n$. Since the map $Aut(F_n)\to GL(n,\Z)$ is surjective, every basis for $\Z^n$ lifts to a basis for $F_n$, and hence every flag of summands of $\Z^n$ lifts to a flag of free factors of $F_n$, i.e. the map $FC_n\to X_n$ is surjective.

Given a basis $\{v_1,\cdots,v_n\}$ for $\Q^n$,  consider the subcomplex of $X_n$ consisting of all flags of subspaces of the form $\langle v_{i_1}\rangle\subset \langle v_{i_1},v_{i_2}\rangle\subset\cdots\subset\langle v_{i_1},\cdots,v_{i_{n-1}}\rangle$.  This subcomplex can be identified with the barycentric subdivision of the boundary of an $(n-1)$-dimensional simplex, so forms an $(n-2)$-dimensional sphere in $X_n$, called an {\it apartment}. 

The construction above applied to  a basis of $F_n$ instead of $\Q^n$ yields an $(n-2)$-dimensional sphere in $FC_n$.
In particular, $H_{n-2}(FC_n)$ is non-trivial.  This sphere maps to an apartment in $X_n$ showing that the induced map $H_{n-2}(FC_n)\to H_{n-2}(X_n)$ is also non-trivial. 

The property of buildings which is missing in $FC_n$ is that given any two maximal simplices there is an apartment which contains both of them.  For example, for $n=3$ and $F_3$ free on $\{x,y,z\}$ there is no ``apartment" which contains both the one-cells corresponding to $\langle x \rangle\subset \langle x,y\rangle$ and 
$\langle yxy\inv\rangle\subset \langle x,y \rangle$, 
since $x$ and $yxy\inv$ do not form part of a basis of $F_3$.

\end{document}